\documentclass[12pt]{article}

\usepackage{algorithm, algpseudocode}
\usepackage{float}
\usepackage{tikz}
\usepackage{capt-of}
\usepackage{color}
\usepackage{fullpage}
\usepackage{amssymb}
\usepackage{amsmath}
\usepackage{amsthm}
\usepackage{amsfonts}
\usepackage{amscd}

\usepackage[colorlinks=true,
linkcolor=brown,
filecolor=brown,
citecolor=brown]{hyperref}

\newcommand{\seqnum}[1]{\href{https://oeis.org/#1}{\rm \underline{#1}}}
\newcommand{\fixp}[2]{\mathrm{Fix}\left(#1,#2\right)}

\begin{document}

\theoremstyle{definition}
\newtheorem{theorem}{Theorem}
\newtheorem{corollary}[theorem]{Corollary}
\newtheorem{lemma}[theorem]{Lemma}
\newtheorem{proposition}[theorem]{Proposition}
\newtheorem{example}[theorem]{Example}
\newtheorem{remark}[theorem]{Remark}

\begin{center}
\vskip 1cm{\LARGE\bf 
Congruences for Dedekind numbers}
\vskip 1cm
\large

Bart\l{}omiej Pawelski\\ Institute of Informatics\\ University of Gda\'nsk\\ 
\href{mailto:bpawelski12@gmail.com}{\tt bpawelski12@gmail.com}

\vskip 1cm
\large
Andrzej Szepietowski\\
Institute of Informatics\\
University of Gda\'nsk\\
\href{mailto:andrzej.szepietowski@ug.edu.pl}{\tt andrzej.szepietowski@ug.edu.pl} \\

\vskip 1cm
\large
B{\l}a{\.z}ej Szepietowski\\
Institute of Mathematics\\
University of Gda\'nsk\\
\href{mailto:blazej.szepietowski@ug.edu.pl}{\tt blazej.szepietowski@ug.edu.pl} \\
\end{center}

\vskip .3in



\begin{abstract}
We study congruence properties of Dedekind numbers $d_n$. Our approach relies on actions of $p$-subgroups of symmetric groups on posets of monotone functions and leads to the congruences $d_{10}\equiv 106\pmod{210}$, $d_{11}\equiv 117\pmod{210}$, and $d_{12}\equiv 12\pmod{30}$. In particular, $d_{11}$ is odd, confirming a claim stated by Yamamoto without a proof.
\end{abstract}
    
\section{Introduction}\label{sec:intro}
Let $B$ denote the set $\{0,1\}$ and $B^n$ the
set of $n$-element sequences of $B$.
A Boolean function with $n$ variables is any function
from $B^n$ into $B$. There is the order relation in $B$ (namely: $0\le 0$, $\;0\le1$, $\;1\le1$)
and the partial order in $B^n$:
for any two elements:
$x=(x_1,\dots,x_n)$, $\;y=(y_1,\dots,y_n)$ in $B^n$,
$x\le y$ if and only if $x_i\le y_i$ for all $1\le i\le n$.
A function $h:B^n\to B$ is monotone if $x\le y$ implies $h(x)\le h(y)$.
Let $D_n$ denote the set of monotone functions with $n$ variables.
We have the partial order in $D_n$ defined by:
$$g\le h\quad\hbox{if and only if}\quad g(x)\le h(x)\quad
\hbox{for all } x\in B^n.$$
Let $d_n$ represent the cardinality of $D_n$, also known as the $n$--th Dedekind number. Dedekind numbers are described by the \textit{On-Line Encyclopedia of Integer Sequences} (OEIS) sequence \seqnum{A000372}. The values $d_0,\dots,d_4$ were found by Dedekind, $d_5$ by Church, $d_6$ by Ward, $d_7$ by Church, and $d_8$ by Wiedemann in 1991~\cite{Wied1991}; the last was confirmed in 2001 by Fidytek, Mostowski, Somla, and Szepietowski~\cite{Fid2001}. In 2023, the value $d_9$ was announced independently by J\"akel~\cite{jakel} and by Van Hirtum et al.~\cite{hirtum}. Dedekind number $d_n$ also corresponds to the number of simple games with $n$ players in minimal winning form, the number of antichains of subsets of an $n$ set, the number of Sperner families and the cardinality of a free distributive lattice on $n$ generators \cite{oeis}.

In this paper we study divisibility properties of Dedekind numbers. In \cite{yamamoto} Yamamoto proved:
\begin{theorem}\label{Y}\cite{yamamoto}
    If $n$ is even, then $d_n$ is also even.
\end{theorem}
Pawelski and Szepietowski \cite{divisibility} proved that:
\[
d_9\equiv 0 \pmod 2,\quad
d_9\equiv 0 \pmod 3,\quad
d_9\equiv 1 \pmod 5,\quad
d_9\equiv 6 \pmod 7.
\]
These results were announced in arXiv before the value $d_9$ was published by J\"akel and Hirtum et al. \cite{jakel, hirtum}. We note that the congruences for $d_9$ obtained by the same method in \cite{divisibility} were later confirmed by the computed value of $d_9$.

Our approach is based on group actions and the orbit decomposition. If a $p$-group acts on a finite set, then the number of fixed points is congruent modulo $p$ to the cardinality of the set (Lemma \ref{p_group_action}). We apply this principle to actions of $p$-subgroups of $S_n$ on posets of monotone functions induced by permuting the variables. For cyclic actions the fixed points are counted using methods developed by the first two authors in \cite{pawelski9,divisibility,szepietowski}.
For non-cyclic actions, we identify the fixed points with monotone functions on an associated quotient poset (Lemma \ref{Ginv}), which is considerably smaller, and more importantly, has much smaller anti-chains than the original domain, leading to a dramatic reduction in the number of monotone functions.


We obtain the following results:
\begin{itemize}
\item $d_{10}\equiv 1\pmod p$ for $p=3,5,7$ (Theorem \ref{thm:d10});  
\item $d_{11}\equiv 1\pmod 2$, $d_{11}\equiv 0\pmod 3$, $d_{11}\equiv 2 \pmod 5$, and $d_{11}\equiv 5\pmod 7$ (Theorem~\ref{thm:d11} and Theorem~\ref{d11mod5}); 
\item $d_{12}\equiv 0\pmod 3$ and $d_{12}\equiv 2\pmod 5$ (Theorem \ref{thm:d12} and Theorem~\ref{d12mod5}). 
\end{itemize}
Combining these results with Theorem \ref{Y}, and using the Chinese Reminder Theorem, we obtain the following congruences.
\[d_{10}\equiv 106\pmod{210},\qquad d_{11}\equiv 117\pmod{210},\qquad d_{12}\equiv 12\pmod{30}.\]
In particular, $d_{11}$ is odd, which confirms the claim given without proof by Yamamoto in \cite{yamamoto}.

The congruences obtained in this paper provide useful consistency checks for future computations of the exact values of $d_{10}$, $d_{11}$ and $d_{12}$.

The methods developed in this paper are not restricted to the cases considered. The same approach can be used to obtain congruences for $d_n$ modulo primes $p\leq n$ for arbitrary $n$. The restriction to $n\leq 12$ and small primes in the present work is due only to computational limitations.

\section{Posets}
A poset (partially ordered set) $(S,\le)$ consists of a set $S$ together with a binary relation (partial order) $\le$ which is reflexive, transitive and antisymmetric.

For two posets $(S,\le)$ and $(T,\le)$,  we denote by $S\times T$ their Cartesian product equipped with with the product order defined by: 
\[(a,b)\le(c,d)  \iff a\le c \text{ and }b\le d.\]
A function $f:S \to T$ is called monotone, if for all $x,y\in S$,  $x \le y$ implies $f(x)\le f(y)$. 
By $T^S$ we denote  the poset of monotone functions from $S$ to $T$ with the partial order  defined by:
$$f\le g\quad\hbox{ if and only if}\quad f(x)\le g(x)\hbox{ for all }x\in S.$$
If there exists a monotone bijection $f:S \to T$, then we say that $S$ and $T$ are isomorphic and write $S\cong T$.

In this paper, we consider the following posets:
\begin{itemize}
\item Path (or chain) $P_m=\{0,1,\dots,m-1\}$ ordered by the usual order on the natural numbers. 
\item The $n$-fold Cartesian product $P_m^n$ equipped with the product order:
\[(x_1,\dots,x_n)\le(y_1,\dots,y_n) \iff x_i\le y_i \text{ for all } i=1,2\dots,n.\]
\item $B=P_2=\{0,1\}$.
\item $B^n$ -- the $n$-dimensional Boolean lattice.
\item $D_n=B^{B^n}$ -- the poset of all monotone Boolean functions of $n$ variables.
\end{itemize}
Note that $D_0=B$ and $D_1\cong P_3$.

In this paper we use the following two lemmas:
\begin{lemma}\label{L1}\cite{szepietowski}
  For three posets $R,S,T$:
The poset $R^{S\times T}$ is isomorphic to $(R^S)^T$ and to $(R^T)^S$.
\end{lemma}

\begin{lemma}\label{L2}\cite{szepietowski}
The poset $D_{n+k}$ is isomorphic to the poset
$D_n^{B^k}$---the poset of monotone functions from
$B^k$ to $D_n$.
\end{lemma}

\section{Group actions and quotient posets}
If a group $G$ acts on a set $X$, we write  
\[\fixp{G}{X}=\{x\in X \colon  gx=x\ \forall g\in G\}\] 
for the set of fixed points.
The following elementary lemma is the principal tool used throughout this paper. Every divisibility result is obtained by applying this lemma to appropriately chosen group action.
\begin{lemma}\label{p_group_action}
Let $p$ be a prime and let $G$ be a finite $p$-group acting on a finite set $X$. Then
\[|\fixp{G}{X}|\equiv|X|\pmod p.\]
\end{lemma}
\begin{proof}
By the Orbit--Stabilizer Theorem, the cardinality of every orbit divides
$|G|$. Since $|G|$ is a power of $p$, every nontrivial orbit has cardinality
divisible by $p$. As the orbits form a partition of $X$, it follows that
$|X|$ is congruent modulo $p$ to the number of singleton orbits, which are precisely  the fixed points of the action.
\end{proof}

Let $S_n$ denote the symmetric group on $\{1,\dots,n\}$ and let $P=P_m$ be a finite path. We consider the natural action of $S_n$ on $P^n$  by permutation
of coordinates:
\[
\sigma(x_1,\ldots,x_n)
=
(x_{\sigma^{-1}(1)},\ldots,x_{\sigma^{-1}(n)}).
\]
The poset $P^n$ is graded with rank function
\[
\rho(x_1,\ldots,x_n)=x_1+\cdots+x_n.
\]
Since permutations of coordinates preserve $\rho$, the action of every subgroup $G\le S_n$ is 
rank-preserving, that is $\rho(\sigma x)=\rho(x)$ for every $\sigma\in G$ and $x\in P^n$. 

For  a subgroup $G\leq S_n$ we define the {quotient poset} $P^n/G$ as follows. The elements of $P^n/G$ 
are the orbits of the action of $G$ on $P^n$. The order relation is induced
by the order on $P^n$:
\[
Gx\leq Gy\ 
\Longleftrightarrow\ 
\exists\, x'\in Gx,\,y'\in Gy \text{ such that }x'\leq y'.
\]
Reflexivity and transitivity follow immediately from the corresponding properties of the order on $P^n$. The antisymmetry follows from the fact that 
$Gx\le Gy$ implies $\rho(x)\le\rho(y)$. Since comparable elements of $P^n$ of the same rank must be equal, two distinct orbits cannot be comparable in both directions.

The following are examples of quotient posets of $B^n$ by an action of a subgroup of $S_n$.
\begin{example}\label{quotient_Bn_Sn}
For the full symmetric group $S_n$,  the quotient poset $B^n/S_n$ is isomorphic to the chain $P_{n+1}$.
Indeed, the rank function induces a bijection between the elements of  $B^n/S_n$ and the set $\{0,1,\dots,n\}$, and this bijection is in fact an isomorphism of posets.
\end{example}
\begin{example}\label{quotient_B3_Z3}
For $C_3=\langle (1\,2\,3)\rangle\leq S_3$, we have $B^3/C_3\cong P_4$. Indeed, the orbits are again determined by the rank, since the action is transitive on each level set of the rank function.
\end{example}
\begin{example}\label{quotient_B3k_Z3k}
Let  \[H=\langle (1\,2\,3),(4\,5\,6),\dots,(3k-2\ 3k-1\ 3k)\rangle\leq S_{3k}.\] Then $H\cong C_3^k$, and the action of $H$ decomposes into independent cyclic actions on the $k$ blocks of three coordinates. Consequently,
\[B^{3k}/H\cong \left(B^3/C_3\right)^k\cong P_4^k.\] 
\end{example}

\section{$G$-invariant monotone functions}
Let $P=P_m$ be a finite path, $T$ an  arbitrary poset, and $G\leq S_n$.
Then $G$ acts on the poset of monotone maps $T^{P^n}$ by
precomposition:
\[
(\sigma f)(x)=f(\sigma^{-1}x),
\qquad \sigma\in G,\ f\in T^{P^n}.
\]
The fixed points of this action are precisely the $G$-invariant
monotone maps:
\[
\fixp{G}{T^{P^n}}=\{f\in T^{P^n} \colon  f(\sigma x)=f(x)\ \forall \sigma\in G, x\in P^n\}.
\]
The following lemma identifies the poset of $G$-invariant monotone maps on $P^n$ with the poset of 
monotone maps on the quotient poset $P^n/G$.

\begin{lemma}\label{Ginv}
Let $P=P_m$ be a finite path, $T$ an  arbitrary poset, and $G\leq S_n$. Then there exists an isomorphism of posets
\[
\fixp{G}{T^{P^n}}\cong T^{P^n/G}.
\]
\end{lemma}
\begin{proof}
For  $f\in\fixp{G}{T^{P^n}}$ we define $\widetilde{f}\in T^{P^n/G}$ by $\widetilde{f}(Gx)=f(x)$. Since $f$ is constant on each orbit, $\widetilde{f}$ is well defined. Moreover, $\widetilde{f}$ is monotone. Indeed, if $Gx\le Gy$, then there exist $x'\in Gx$, $y'\in Gy$ such that $x'\le y'$, and
\[\widetilde{f}(Gx)=f(x)=f(x')\le f(y')=f(y)=\widetilde{f}(Gy).\] 
Define \[\Phi\colon\fixp{G}{T^{P^n}}\to T^{P^n/G},\qquad \Phi(f)=\widetilde{f}.\] The inverse map is given by \[\Phi^{-1}(\widetilde{f})(x)=\widetilde{f}(Gx).\]
Finally, $\Phi$ preserves the pointwise order, since \[f_1\le f_2\Longleftrightarrow \widetilde{f_1}\le\widetilde{f_2}.\] Therefore, $\Phi$ is an isomorphism of posets.
\end{proof}

\begin{lemma}\label{p_quotient}
Let $P=P_m$ be a finite path,  and $G$ a $p$-subgroup of $S_n$. Then, for arbitrary poset $T$,
\[\left|T^{P^n}\right|\equiv\left| T^{P^n/G}\right| \pmod p.\]
\end{lemma}
\begin{proof}
By applying Lemma \ref{p_group_action} to the action of $G$ on $T^{P^n}$ we have 
\[\left|T^{P^n}\right|\equiv\left| \fixp{G}{T^{P^n}}\right| \pmod p,\]
and by Lemma \ref{Ginv}, \[\left|\fixp{G}{T^{P^n}}\right|=\left| T^{P^n/G}\right|.\]
\end{proof}
\section{Congruences from cyclic $p$-group actions}
Recall, that the $n$-th Dedekind number is $d_n=|D_n|$, where $D_n=B^{B^n}$.
In this section we compute $d_{10}\bmod p$ for $p=3,5,7$ and $d_{11}\bmod p$ for $p=2,3,7$ by applying Lemma \ref{p_group_action} to cyclic $p$-group actions on $D_n$, $n\in\{10,11\}$.

For the rest of this section let us introduce the following notation. For $\sigma\in S_n$ we denote by $\phi_n(\sigma)$ the cardinality of the poset of $\sigma$-invariant monotone functions from $B^n$ to $B$, that is
\[\phi_n(\sigma)=|\fixp{\langle\sigma\rangle}{D_n}|=|\{f\in D_n\colon f(\sigma x)=f(x)\ \forall x\in B^n\}|.\]
If the order of $\sigma$ is a power of a prime $p$, then by Lemma \ref{p_group_action}, 
\begin{equation}\label{eq:dmodp}
\phi_n(\sigma)\equiv d_n\pmod p.
\end{equation}

\begin{theorem}\label{thm:d10}
The number $d_{10}$ satisfies the congruence
\[d_{10}\equiv 1\pmod p\quad\text{for }p=3,5,7.\]
 \end{theorem}
 \begin{proof}
 We computed the values of of $\phi_{10}(\sigma)$, where $\sigma\in S_{10}$ is a cycle of length $5$, $7$ and $9$. By applying \eqref{eq:dmodp} we obtain:
 \begin{align*}
&d_{10}\equiv\phi_{10}((1\,2\,3\,4\,5))=16\,578\,128\,790\,019\,768\,076\equiv 1\pmod 5.\\
&d_{10}\equiv\phi_{10}((1\,2\,3\,4\,5\,6\,7))= 1\,043\,977\,055\,227\equiv 1\pmod 7.\\
&d_{10}\equiv\phi_{10}((1\,2\,3\,4\,5\,6\,7\,8\,9))=1\,702\,659\,955\equiv 1\pmod 3.
\end{align*}
\end{proof}
The residues of $d_{10}$  modulo $3$ and $5$ will be confirmed by a different method in the next section; see Remarks \ref{d_10mod3} and \ref{d_10mod5} respectively.
\begin{theorem}\label{thm:d11}
The number $d_{11}$ satisfies the congruences
\[d_{11}\equiv 1\pmod 2,\quad d_{11}\equiv 0\pmod 3,\quad d_{11}\equiv 5\pmod 7.\]
 \end{theorem}
 \begin{proof}
 We computed the values of of $\phi_{11}(\sigma)$, where $\sigma\in S_{11}$ is a cycle of length $7$, $8$ and $9$. By applying \eqref{eq:dmodp} we obtain:
 \begin{align*}
&d_{11}\equiv\phi_{11}((1\,2\,3\,4\,5\,6\,7))=18\,035\,464\,634\,257\,554\,417\,288\equiv 5\pmod 7.\\
&d_{11}\equiv\phi_{11}((1\,2\,3\,4\,5\,6\,7\,8))= 345\,413\,653\,441\,343\,650\,195\equiv 1\pmod 2.\\
&d_{11}\equiv\phi_{11}((1\,2\,3\,4\,5\,6\,7\,8\,9))=126\,088\,909\,171\,123\,806\equiv 0\pmod 3.
\end{align*}
\end{proof}
The numbers $\phi_{11}\big((1\,2\,3\,4\,5)\big)$ and  $\phi_{11}\big((1\,2\,3\,4\,5)(6\,7\,8\,9\,10)\big)$
lie beyond our current computational reach. However, in the next section we compute $d_{11}\bmod 5$ by a different approach, using a non-cyclic group of order $25$ -- see Theorem \ref{d11mod5}.
\section{Congruences from non-cyclic $p$-group actions}
In this section we use non-cyclic group actions and the notion of a quotient poset to study the divisibility of Dedekind numbers.

 \begin{lemma}\label{dn+3kmod3}
For $n\geq 0$ and $k\geq 0$ we have
\[d_{n+3k}\equiv |D_n^{P_4^k}|\pmod{3}.\] 
\end{lemma}
\begin{proof}
By Lemma \ref{L2},\[D_{n+3k}\cong D_n^{B^{3k}}.\]
Let $H\leq S_{3k}$ be the group from Example \ref{quotient_B3k_Z3k}, of order $3^k$. We have \[B^{3k}/H\cong P_4^k,\] and by Lemma \ref{p_quotient},
\[\left|D_n^{B^{3k}}\right|\equiv\left|D_n^{P_4^k}\right|\pmod 3.\]
\end{proof}

\begin{theorem}\label{thm:d12}
The number $d_{12}$ is divisible by $3$.
\end{theorem}
\begin{proof}
 By Lemma \ref{dn+3kmod3},
 \[d_{12}\equiv \left|B^{P_4^4}\right|\pmod 3.\]
 Since $B^{P_4}\cong P_5$, Lemma \ref{L1} gives
 \[B^{P_4^4}\cong P_5^{P_4^3}.\]
 Consider the action of $C_3=\langle(1\,2\,3)\rangle\le S_3$ on $P_4^3$ by permutation of coordinates and let
  \[Q=P_4^3/C_3\] be the quotient poset of $P_4^3$ by this action. By Lemma \ref{p_quotient},
  \[d_{12}\equiv \left| P_5^{P_4^3}\right|\equiv\left|P_5^Q\right|\pmod 3. \]
    Note that $Q$ has $24$ elements. Indeed, there are $4$ fixed points of the action and $20$ orbits of length $3$. 
 A computer calculation gives
 \[\left|P_5^Q\right|=3\,070\,632,\]
 hence $d_{12}\equiv 0\pmod 3$.
\end{proof}

\begin{remark}\label{d_10mod3}
The method of proof of Theorem \ref{thm:d12} can also be used to verify the computation of $d_{10}\bmod 3$ from Theorem \ref{thm:d10}.
By Lemma \ref{dn+3kmod3} we have
\[d_{10}\equiv\left|D_1^{P_4^3}\right|\pmod 3.\]
Passing to the quotient poset $Q$ from the proof of Theorem \ref{thm:d12}, by Lemma \ref{p_quotient}  we obtain
\[\left|D_1^{P_4^3}\right|\equiv\left|D_1^Q\right|\pmod 3.\]
 A computer calculation gives
 \[\left|D_1^Q\right|=\left|P_3^Q\right|=6526,\]
 hence $d_{10}\equiv 1\pmod 3$.    
\end{remark}
\begin{remark} 
 The poset $Q=P_4^3/C_3$ appearing in proof of Theorem \ref{thm:d12} can equivalently be realized  as the quotient of $B^9$ by the action of the  (non-abelian)  Sylow $3$-subgroup of $S_9$, of order $3^4$, generated by $(1\,2\,3)$, $(4\,5\,6)$, $(7\,8\,9)$ and $(1\,4\,7)(2\,5\,8)(3\,6\,9)$. 
\end{remark}

\begin{lemma}\label{dn+5kmod5}
Let $R$ be the quotient poset of $B^5$ by the action of $C_5=\langle(1\,2\,3\,4\,5)\rangle\le S_5$ by permutation of coordinates.
For $n\geq 0$ and $k\geq 1$ we have
\[d_{n+5k}\equiv |D_n^{R^k}|\pmod{5}.\] 
\end{lemma}
\begin{proof}
By Lemma \ref{L2},\[D_{n+5k}\cong D_n^{B^{5k}}.\]
Let $G\leq S_{5k}$ be the group of order $5^k$ generated by $k$ disjoint $5$-cycles:
\[G=\langle(1\,2\,3\,4\,5), \dots,(5k-4\ 5k-3\ 5k-2\ 5k-1\ 5k)\rangle.\]
We have \[B^{5k}/G\cong\left(B^5/C_5\right)^k\cong R^k,\] and by Lemma \ref{p_quotient},
\[ \left|D_n^{B^{5k}}\right|\equiv\left|D_n^{R^{k}}\right| \pmod{5}.\]
\end{proof}

\begin{theorem}\label{d11mod5}
The number $d_{11}$ satisfies the congruence \[d_{11}\equiv 2\pmod 5.\] 
\end{theorem}
\begin{proof}
By Lemma \ref{dn+5kmod5},
\[d_{11}\equiv|D_1^{R\times R}|=|P_3^{R\times R}|\pmod 5.\]  
Observe that $R$ has $8$ elements.
A computer calculation gives
\[|P_3^{R\times R}|=825\,398\,172,\]
hence, $d_{11}\equiv 2\pmod 5.$
\end{proof}

\begin{theorem}\label{d12mod5}
The number $d_{12}$ satisfies the congruence \[d_{12}\equiv 2\pmod 5.\] 
\end{theorem}
\begin{proof}
By Lemma \ref{dn+5kmod5} applied with $n=2$ and $k=2$,
\[d_{12}\equiv\left|D_2^{R\times R}\right|\pmod 5.\]  
Observe that $D_2$ has $6$ elements and $R\times R$ has $64$ elements.
A computer calculation gives
\[\left|D_2^{R\times R}\right|=37\,570\,480\,680\,918\,252,\]
hence, $d_{12}\equiv 2\pmod 5.$
\end{proof}

\begin{remark}\label{d_10mod5}
Lemma \ref{dn+5kmod5}  can also be used to verify the computation of $d_{10}\bmod 5$ from Theorem \ref{thm:d10}.
 A computer calculation gives
 \[\left|B^{R\times R}\right|=71716,\]
 hence, by Lemma \ref{dn+5kmod5},
\[d_{10}\equiv|B^{R\times R}|\equiv 1\pmod 5.\]      
\end{remark}

\section{Computations}\label{sec:computations}
All the numerical values stated in this paper can be recomputed with a single Java program, available in the repository \cite{congruences_repo}. The program has no external dependencies, verifies every poset it constructs, and compares each computed number with the value stated here. The computations are parallelised; a complete run took about $28$ minutes on a laptop with an AMD Ryzen~7 8840U processor.

\section*{AI Disclosure}
We used Claude Opus 5 to assist in optimizing the Java programs used for the computations reported in this paper. The programs are available in the repository \cite{congruences_repo}. The authors reviewed the AI-assisted modifications and verified the correctness of the programs. The authors retain full responsibility for the mathematical content of the paper.

    \end{document}